\documentclass[a4paper,10pt]{article}

\usepackage{amsmath, amsfonts, amssymb, latexsym, mathrsfs, graphicx}

\usepackage[all]{xy}
\newdir{ ->}{{}*!/-5ex/@{->}}

\newcommand{\mmap}[3]{\ar@/#2/[#1]|*=0{\rotatebox{#3}{$\scriptsize |$}}} %\mmap{Direction}{CurvatureNoneIfEmpty}{Rotation90IsInvariant}%
\newcommand{\picar}[3]{\ar@{}[#1]|*{\rotatebox{#2}{$\scriptsize #3$}}} %\picar{Direction}{Rotation90IsInvariant}{Symbol}%

\newcommand{\sa}{{\mathscr{A}}}

\newcommand{\ssc}{{\mathscr{C}}}
\newcommand{\oop}{\operatorname{op}}
\newcommand{\op}{^{\oop}}
\newcommand{\set}{\mathbf{Set}}
\newcommand{\vect}{\mathbf{Vect}_{k}}
\newcommand{\enc}{\mathbf{End}^{\vee}U}

\begin{document}

\begin{center}
{\underline{{\Large Note On Endomorphism Algebras Of Separable}}}
\begin{center}
{\underline{{\Large Monoidal\phantom{p}Functors.}}}
\end{center}

{{\small Brian J. DAY and Craig A. PASTRO}}

{{\small April 26, 2009.}}
\end{center}

\noindent \emph{\small \underline{Abstract:} We recall the Tannaka construction in [2] for certain types of split monoidal functor into $\vect$, and remove the ``compactness'' restriction on the domain.}
\newline

\begin{center}
\begin{tabular}{c}\hline
\phantom{a}\phantom{a}\phantom{a}\phantom{a}\phantom{a}\phantom{a}\phantom{a}\phantom{a}\phantom{a}\phantom{a}\phantom{a}\phantom{a}\phantom{a}\phantom{a}\phantom{a}\phantom{a}
\end{tabular}
\end{center}

In $\set$ or $\vect$, a ``VN-core'' is an algebra $(A,\mu,\eta)$ and a coalgebra $(A,\delta,\epsilon)$ satisfying the axiom
$$\delta\mu = (\mu\otimes\mu)(1\otimes c\otimes 1)(\delta\otimes\delta),$$
where $c$ is the symmetry for $\otimes$, with a map $S:A\to A$ such that
$$\mu_{3}(1\otimes S\otimes 1)\delta_{3} = 1: A\to A\, .$$
The core is called \underline{unital} if it satisfies the stronger axiom
$$1\otimes\eta = (1\otimes\mu)(1\otimes S \otimes 1)\delta_{3}\, .$$
In $\set$, unital VN-cores are precisily groups, while in $\vect$ they contain the Hopf $k$-algebras.

The following theory is $k$-linear for $k$ a field of characteristic $0$, and we shall assume some familiarity with [2] Theorem 3.1 which we now generalize.
The aim of this procedure is to remove the construction of the unital VN-core $\enc$ in [2] from the ``compact'' setting by postulating the existence of an \underline{abstract} duality functor
$$(-)^{\ast}: \sa\op\to\sa$$
on the generating category $\sa\subset\ssc$ with an associated natural isomorphism 
$$u: U(A^{\ast}) \cong U(A)^{\ast}$$
and a transformation $e:A^{\ast}\otimes A\to I$ satisfying the two axioms
\newline

\noindent\underline{$(e,r,r_{0})$:}
\begin{displaymath}
\xymatrix{
k \ar[rr]^{r_{0}} && UI \\
U(A)^{\ast}\otimes UA \ar[u]^{e} && U(A^{\ast}\otimes A) \ar[u]_{Ue} \\
&U(A^{\ast})\otimes UA \ar[ur]_{r} \ar[ul]^{u\otimes 1}&
}
\end{displaymath}

\noindent commutes, and

\noindent\underline{$(e,i,i_{0})$:}
\begin{displaymath}
\xymatrix{
UI \ar[rr]^{i_{0}} && k \\
U(A^{\ast}\otimes A) \ar[u]^{Ue} \ar[dr]_{i} && U(A)^{\ast}\otimes UA \ar[u]_{e} \\
&U(A^{\ast})\otimes UA \ar[ur]_{u\otimes 1} &
}
\end{displaymath}

\noindent commutes.

We emphasize here that it is not necessary for $A^{\ast}$ to be the actual $\otimes$-dual of $A$ in the monoidal category $\ssc$.
For example, $\ssc$ could be a $\star$-autonomous monoidal category [1], or just a monoidal closed category [3], and so forth. 

The two axioms above now replace the cumbersome ``$U$-trace'' condition in [2] $\S$ 3.
The region  labeled (3) in [2] $\S$ 3 now becomes:

\begin{displaymath}
\def\objectstyle{\scriptstyle} 
\def\labelstyle{\scriptstyle}
\xymatrix@C-0.4cm{
&& (U(A)^{*}\otimes UA)^{*} \otimes U(A)^{*}\otimes UA \ar[rrd]^-{\phantom{a}\phantom{a}\phantom{a}\phantom{a}\phantom{a}(u\otimes 1)^{*}\otimes (u^{-1}\otimes 1)} && \\
k^{*} \otimes U(A)^{*} \otimes UA \ar[rrd]^{\phantom{a}\phantom{a}1\otimes u^{-1}\otimes 1}_<{\phantom{.}}="p1"  \ar[urr]^-{e^{*}\otimes 1 \otimes 1\phantom{a}\phantom{a}} &&&& \ar[dd]^{i^{*}\otimes r} (U(A^{*})\otimes UA)^{*} \otimes U(A^{*})\otimes UA \\
&& k^{*} \otimes U(A^{*}) \otimes UA \ar[dd]^{i^{*}_{0}\otimes r}  \ar@{}[urr]|*=0-{(e,i,i_{0})} && \\
k^{*} \otimes UA \otimes U(A)^{*} \ar[uu]^{1\otimes c} &&&& U(A^{*}\otimes A)^{*} \otimes U(A^{*}\otimes A) \ar[ddd]^{copr_{B =A^{*}\otimes A}} \\
&& \ar@{}[uul]|*=0-{(e,r,r_{0})\phantom{aaaa}} U(I)^{*}\otimes U(A^{*}\otimes A) \ar[rru]_{\phantom{a}\phantom{a}U(e)^{*}\otimes 1} \ar[dd]^{1\otimes U(e)} && \\
& k^{\ast}\otimes k \ar[rd]^{i^{*}_{0}\otimes r_{0}} \ar@{}="p2" \ar"p1";"p2"^{1\otimes e}="mid" &&&\\
k^{*}\otimes k \ar[ru]_{1} \ar@{}="start" \ar@{}"start";"mid"|*=0{(trace)}  \ar[uuu]^{1\otimes (\textrm{dim} UA)^{-1}) .n}  \ar[rr]_-{i_{0}^{*}\otimes r_{0}} && U(I)^{*}\otimes UI  \ar[rr]_-{copr_{B= I}} && \int^{B} U(B)^{*}\otimes UB \ar@{}[uull]_-{nat.}
}
\end{displaymath}

\noindent which eventually commutes by commutativity of

\begin{displaymath}
\xymatrix{
UA\otimes U(A)^{*} \ar[rr]^-{c} && U(A)^{*}\otimes UA \ar[d]^-{e} \\
k \ar[u]^-{n} \ar[rr]_-{(\textrm{dim }UA).1} \ar@{}[rru]|*=0-{(trace)} && k
}
\end{displaymath}
in $\vect$.

\pagebreak

Thus we essentially replace the map $S$ in [2] Theorem 3.1 by the map defined on $\enc$ by the new family:
$$S_{A} = (\textrm{dim }UA)^{-1} . \sigma_{A}$$
(i.e. remove the factor $(\textrm{dim }UI)$ and assume $\textrm{dim }UA \neq 0$ for all $A\in\sa$).
The new structure $(\enc,\mu,\eta,\delta,\epsilon,S)$ is then a unital VN-core in $\vect$.
Note that we no longer require $\ssc$ and $U$ to be \underline{braided} structures in this context.
\newline
\newline

\begin{center}
{\underline{{\Large References.}}}
\end{center}

\begin{center}
\begin{tabular}{l l}
\textrm{[1]} & M. Barr, ``$\star$-Autonomous Categories", \\
 & Lecture Notes in Mathematics, (Springer) 752 (1979).\\
 & \\
\textrm{[2]} & B. J. Day and C. A. Pastro, ``On Endomorphism Algebras of Separable Monoidal Functors'', \\
 & Theory Appl. Categories, Vol. 22, No. 4 (2009), Pg. 77-96. (Electronic).  \\
 & \\
\textrm{[3]} & S. Eilenberg  and G. M. Kelly, ``Closed Categories'', \\ 
 & Proc. Conf. on Categorical Algebra, La Jolla 1965, Springer-Verlag (1966) Pg. 421-562. \\
 & \\
 & \\
 & \\
 & \\
\end{tabular}

\small{Mathematics Dept., Faculty of Science, Macquarie University, NSW 2109, Australia.}
\newline

Any replies are welcome through Tom Booker (thomas.booker@students.mq.edu.au), who kindly typed the manuscript.

\end{center}

\end{document}